\documentclass[a4paper,12pt]{article}
\usepackage{amsfonts}
\usepackage{enumerate}

\def\Empty{\phi}
\def\N{\mathbb N}
\def\R{\mathbb R}
\def\Rn{\R^n}
\def\Rp{\R^p}

\def\n{\mathcal N}
\renewcommand{\H}{\Rn\x\Rp}

\def\set#1{\left\{#1\right\}}
\def\Ia{\set{1,...,\alpha}}
\def\ip#1{\left<#1\right>} 
\def\ipi#1{\left<#1\right>_i}
\def\norm#1{\left\Vert#1\right\Vert}
\def\normi#1{\norm{#1}_i}
\def\normw#1{\norm{#1}_{\infty,\gamma}}
\def\normmax#1{\norm{#1}_{\infty}}

\def\Lim{\lim\limits}
\def\Prod{\prod\limits}
\def\Sum{\sum\limits}
\def\Min{\min\limits}
\def\Max{\max\limits}
\def\Inf{\inf\limits}
\def\Sup{\sup\limits}
\def\x{\times}

\def\proof{\noindent{\bf{Proof. }}}
\def\endproof{\hfill{\rule{0.5em}{0.5em}}}
\def\keywords#1\par{\\ \newline \textbf{Key words:} #1}

\newtheorem{theoreme}{Theorem}
\newtheorem{corollaire}[theoreme]{Corollary}

\newtheorem{remarque}{Remark}

\begin{document}
\title{\uppercase{Parallel synchronous algorithm for nonlinear fixed point
problems}}
\date{}

\author{
Ahmed ADDOU \\
and\\
Abdenasser BENAHMED}

\maketitle
\begin{abstract}
\begin{footnotesize}
We give in this paper a convergence result concerning parallel
synchronous algorithm for nonlinear fixed point problems with
respect to the euclidian norm in $\Rn$. We then apply this result
to some problems related to convex analysis like minimization of
functionals, calculus of saddle point, convex programming...
\end{footnotesize}
\end{abstract}

\section{Introduction.}\label{int}
This study is motivated by the paper of Bahi[\ref{bahi99}] where
he has given a convergence result concerning parallel synchronous
algorithm for linear fixed point problems using nonexpansive
linear mappings with respect to a weighted maximum norm. Our goal
is to extend this result to a nonlinear fixed point problems,
\begin{equation}\label{pointfixe}
F(x^*)=x^*
\end{equation}
with respect to the euclidian norm,
where $F:\Rn\to\Rn$ is a nonlinear operator. \\
Section \ref{paa} is devoted to a brief description of
asynchronous parallel algorithm. In section \ref{cag} we prove the
main result concerning the convergence of the general algorithm in
the synchronous case to a fixed point of a nonlinear operator from
$\Rn$ to $\Rn$. A particular case of this algorithm (Algorithm of
Jacobi) is applied in section \ref{app} to the operator
$F=(I+T)^{-1}$ which is called the proximal mapping associated
with the maximal monotone operator $T$ (see
Rockafellar[\ref{fellar76a}]).
\section{Preliminaries on asynchronous algorithms.}\label{paa}
Asynchronous algorithms are used in the parallel treatment of
problems taking in consideration the interaction of several
processors. Write $\Rn$ as the product $\Prod_{i=1}^\alpha
\R^{n_i}$, where $\alpha\in\N-\{0\}$ and $n=\Sum_{i=1}^{\alpha
}{n_i}$. All vectors $x\in\Rn$ considered in this study are
splitted in the form $x=(x_1,...,x_\alpha)$ where
$x_i\in\R^{n_i}$. Let $\R^{n_i}$ be equipped with the inner
product $\ipi{.,.}$ and the associated norm
$\normi{.}=\ipi{.,.}^{1/2}$. $\Rn$ will be equipped with the inner
product $\ip{x,y}=\Sum_{i=1}^{\alpha }\ipi{x_i,y_i}$ where
$x,y\in\Rn$ and the associated norm
$\norm{x}=\ip{x,x}^{1/2}=(\Sum_{i=1}^{\alpha}\normi{x_i}^2)^{1/2}$.
It will be equipped also with the maximum norm defined by,
$$
\normmax{x}=\Max_{1\leq i\leq \alpha}\normi{x_i}
$$
Define :\\
$J=\set{J(p)}_{p\in\N}$ a sequence of non empty sub sets of $\Ia$ and \\
$S=\set{(s_1(p),...,s_\alpha(p))}_{p\in\N}$ a sequence of
$\N^\alpha$ such that,
\begin{enumerate}[$\bullet$]
\item $\forall i \in \Ia$, the subset $\{p\in\N,i\in J(p)\}$ is infinite.
\item $\forall i \in \Ia,\forall p\in \N,s_i(p)\le p.$
\item $\forall i \in \Ia, \Lim_{p\to\infty } s_i(p)=\infty.$
\end{enumerate}
Consider an operator $F=(F_1,...,F_\alpha):\Rn\to\Rn$ and define
the asynchronous algorithm associated with $F$ by,
\begin{equation}\label{algorithme}
\left\{
\begin{array}{l}
x^0=(x_1^0,...,x_{\alpha }^0)\in \Rn \\
x_i^{p+1}=\left\{
\begin{array}{lll}
x_i^p  &  if  &  i\notin J(p)\\
F_i(x_1^{s_1(p)},...,x_{\alpha}^{s_{\alpha}(p)})  &  if  & i\in
J(p)
\end{array}
\right. \\
i=1,...,\alpha\\
p=0,1,..
\end{array}
\right.
\end{equation}
(see Bahi and al.[\ref{bahi97}], El Tarazi[\ref{eltarazi82}]). It
will be denoted by $(F,x^0,J,S)$. This algorithm describes the
behaviour of iterative process executed asynchronously on a
parallel computer with $\alpha$ processors. At each iteration
$p+1$, the $i^{th}$ processor computes $x_i^{p+1}$ by using
(\ref{algorithme}) (Bahi[\ref{bahi98}]).\\
$J(p)$ is the subset of the indexes of the components updated at the $p^{th}$ step.\\
$p-s_i(p)$ is the delay due to the $i^{th}$ processor when it computes
the $i^{th}$ block at the $p^{th}$ iteration.\\
If we take $s_i(p)=p \ \forall i\in \Ia$, then (\ref{algorithme})
describes synchronous algorithm (without delay). During each
iteration, every processor executes a number of computations that
depend on the results of the computations of other processors in
the previous iteration. Within an iteration, each processor does
not interact with other processors, all interactions takes place
at the end of iterations
(Bahi[\ref{bahi99}]).\\
If we take
\[
\left\{
\begin{array}{ll}
s_i(p)=p   &  \forall p \in \N,\forall i \in \Ia\\
J(p)=\Ia  &  \forall p \in \N
\end{array}
\right. \\
\]
then (\ref{algorithme}) describes the algorithm of Jacobi.\\
If we take
\[
\left\{
\begin{array}{ll}
s_i(p)=p   &  \forall p \in \N,\forall i \in \Ia\\
J(p)=p+1 \ (mod\ \alpha)   &  \forall p \in \N
\end{array}
\right. \\
\]
then (\ref{algorithme}) describes the algorithm of Gauss-Seidel.\\
For more details about asynchronous algorithms see [\ref{bahi97}],
[\ref{bahi98}], [\ref{bahi99}] and [\ref{eltarazi82}].\\
In the following theorem, Bahi[\ref{bahi99}] has shown the
convergence of the sequence $\{x^p\}$ defined by
(\ref{algorithme}) in the synchronous linear case, $i.e$ $F$ is a
linear operator and $s_i(p)=p, \ \forall p\in\Ia$.
\begin{theoreme}
Consider $\set{T^p}_{p\in\N}$ a sequence of matrices in $\R^{n\x
n}$.
Suppose\\
$(h_{0})$ $\exists$ a subsequence $\{p_{k}\}_{k\in \N}$ such that
$J(p_{k})=\Ia,$ \\
$(h_{1})$\ $\exists\gamma \gg 0\footnote{ $\ \gamma \gg 0$ means
 $\gamma_i>0\ \forall i\in\Ia$ },\forall p\in \N,$\ \
$T^p$\ is nonexpansive\footnote{ \ A matrice $A\in\R^{n\x n}$ is
said to be nonexpansive with respect to the norm $\norm{.}$ if
$\forall x\in \Rn,\ \norm{Ax}\le\ \norm{x}$. $A$ is said to be
paracontracting if $\forall x\in \Rn,\ x\neq Ax \iff\norm{Ax}<
\norm{x}$. } with respect to a weighted maximum norm $\normw{.}$
defined by
\[
x\in \Rn,\ \normw{x}=\max_{1 \le i \le \alpha}{\frac{\normi{x_i}}{\gamma_i}}
\]
$(h_{2})$ $\set{T^p}_{p\in\N}$\ converges to a matrix $Q$ which is
paracontracting with respect to the norm $\normw{.}$.\\
$(h_{3})$\ $\forall p\in \N,\ \n(I-Q)\subseteq \n(I-T^{p})$ ($\n$
denotes the null space).
\\ then
\begin{enumerate}
\item  $\forall x^0\in \Rn$ the sequence $\set{x^p}_{p\in\N}$\ is convergent in $\Rn$
\item  $\Lim_{p\to \infty }x^p=x^*\in \n(I-Q)$
\end{enumerate}
\end{theoreme}
\proof
See Bahi[\ref{bahi99}].
\endproof

\begin{remarque}
The hypothesis $(h_0)$ means that the processors are synchronized
and all the components are infinitely updated at the same
iteration. This subsequence can be chosen by the programmer.
\end{remarque}

\section{Convergence of the general algorithm.}\label{cag}
We establish in this section the convergence of the general
parallel synchronous algorithm to a fixed point of a nonlinear
operator $F:\Rn\to\Rn$ with respect to the euclidian norm defined
in section \ref{paa}. We recall that an operator $F$ from $\Rn$ to
$\Rn$ is said to be nonexpansive with respect to a norm $\norm{.}$
if,
$$
\forall x,y\in\Rn,\ \norm{F(x)-F(y)} \le \norm{x-y}
$$
\begin{theoreme}\label{thprincipal}
Suppose\\
$(h_0)\ \exists$ a subsequence $\{p_{k}\}_{k\in \N}$ such that
$J(p_{k})=\Ia\\
(h_1)\ \exists u\in\Rn,\ F(u)=u \\
(h_2)\ \forall x,y\in\Rn,\ \normmax{F(x)-F(y)} \le \normmax{x-y} \\
(h_3)\ \forall x,y\in\Rn,\ \norm{F(x)-F(y)}^2 \le \ip{F(x)-F(y),x-y}\\
$ Then any parallel synchronous\footnote{ \ In this case
$s_i(p)=p\ \forall i\in \Ia\ \forall p\in\N$.}algorithm defined by
(\ref{algorithme}) associated with the operator $F$ converges to a
fixed point $x^*$ of $F$.
\end{theoreme}
\proof
\begin{enumerate}[(\it i)]
\item\label{etape1} We prove first that the sequence $\set{x^p}_{p\in\N}$
is bounded.\\
$\forall i \in \Ia$ we have,\\
either $i \notin J(p)$, so
\[
\begin{tabular}{lll}
$\normi{x_i^{p+1}-u_i}$
&  =     & $\normi{x_i^p-u_i}$ \\
&  $\le$ & $\normmax{x^p-u}$\\
\end{tabular}\\
\]
or $i \in J(p)$, so
\[
\begin{tabular}{lll}
$\normi{x_i^{p+1}-u_i}$
&  =     & $\normi{F_i(x^p)-F_i(u)}$ \\
&  $\le$ & $\normmax{F(x^p)-F(u)}$\\
&  $\le$ & $\normmax{x^p-u}\ (by\ (h_2))$\\
\end{tabular}\\
\]
so
$$
\forall i\in\Ia,\ \normi{x_i^{p+1}-u_i} \le \normmax{x^p-u}
$$
then
$$
\forall p\in\N,\ \normmax{x^{p+1}-u} \le \normmax{x^p-u}
$$
hence
$$
\forall p\in\N,\ \normmax{x^p-u} \le \normmax{x^0-u}
$$
this proves that the sequence $\set{x^p}_{p\in\N}$ is bounded with
respect the maximum norm and then it's bounded with respect the
euclidian norm .
\item
As the sequence $\set{x^{p_k}}_{k\in\N}$ is bounded
($\set{p_k}_{k\in\N}$ is defined by $(h_0)$), it contains a
subsequence noted also $\set{x^{p_k}}_{k\in\N}$ which is
convergent in $\Rn$ to an $x^*$. We show that $x^*$ is a fixed
point of $F$. For it, we consider the sequence
$\{y^p=x^p-F(x^p)\}_{p\in \N}$ and prove that
$\Lim_{k\to\infty}{y^{p_k}}=0$.
\[
\begin{tabular}{lll}
$\norm{x^{p_k}-u}^2$
&  =  & $\norm{y^{p_k}+F(x^{p_k})-u}^2$\\
&  =  & $\norm{y^{p_k}}^2+\norm{F(x^{p_k})-u}^2+2\ip{F(x^{p_k})-u,y^{p_k}}$\\
\end{tabular}
\]
however
\[
\begin{tabular}{lll}
$\ip{F(x^{p_k})-u,y^{p_k}}$
&  =  & $\ip{F(x^{p_k})-F(u),x^{p_k}-F(x^{p_k})}$\\
&  =  & $\ip{F(x^{p_k})-F(u),[x^{p_k}-F(u)]-[F(x^{p_k})-F(u)]}$\\
&  =  & $\ip{F(x^{p_k})-F(u),x^{p_k}-u}-\norm{F(x^{p_k})-F(u)}^2$\\
&  $\geq$  & $0\ (by\ (h_3))$\\
\end{tabular}
\]
so,
\[
\begin{tabular}{lll}
$\norm{y^{p_k}}^2$
&  $\le$  & $\norm{x^{p_k}-u}^2 - \norm{F(x^{p_k})-u}^2$\\
&  =       & $\norm{x^{p_k}-u}^2 - \norm{x^{p_k+1}-u}^2\ (by\ (h_0))$\\
\end{tabular}
\]
However, in (\textit{\ref{etape1}}) we have shown in particular
that the sequence $\set{\normmax{x^p-u}}_{p\in\N}$ is decreasing
(and it's positive), it's therefore convergent, then the sequence
$\set{\norm{x^p-u}}_{p\in\N}$ is also convergent, so
\[
\begin{tabular}{lll}
$\Lim_{p\to\infty}{\norm{x^p-u}}$
&  =  & $\Lim_{k\to\infty}{\norm{x^{p_k}-u}}$\\
&  =  & $\Lim_{k\to\infty}{\norm{x^{p_k+1}-u}}$\\
&  =  & $\norm{x^*-u}$\\
\end{tabular}
\]
and so
$$
\Lim_{k\to\infty}{\norm{y^{p_k}}}=0
$$
which implies that
$$
x^*-F(x^*)=0
$$
that is $x^*$ is a fixed point of $F$.
\item We prove as in (\textit{\ref{etape1}}) that the sequence
$\set{\normmax{x^p-x^*}}_{p\in\N}$ is convergent, so
$$
\Lim_{p\to\infty}{\normmax{x^p-x^*}}
=\Lim_{k\to\infty}{\normmax{x^{p_k}-x^*}}=0
$$
Which proves that $x^p\to x^*$ with respect to the uniform norm
$\normmax{..}$.
\endproof
\end{enumerate}
\begin{remarque}
We have used the hypothesis $(h_2)$ to prove that the sequence
$\set{x^p}_{p\in\N}$ is bounded. In the case of the parallel
algorithm of Jacobi where $J(p)=\Ia\ \forall p\in\N$, we don't
need this hypothesis, since in this case $x ^{p+1}=F(x^p) \
\forall p\in\N$, and use $(h_3)$ to obtain
$$
\norm{x^{p+1}-u}=\norm{F(x^p)-F(u)} \le  \norm{x^p-u},
$$
hence the corollary,
\end{remarque}

\begin{corollaire}\label{jacobi}
Under the hypotheses $(h_1)$, $(h_3)$ and\\
$(h_0^{'})$ $\forall p\in\N,\ J(p)=\Ia$\\
The parallel Jacobi algorithm defined by
\begin{equation}
\left\{
\begin{array}{l}
x^0=(x_1^0,...,x_{\alpha}^0)\in\Rn \\
x_i^{p+1}=F_i(x_1^p,...,x_{\alpha}^p)\\
i=1,...,\alpha\\
p=1,2...\\
\end{array}
\right.\\
\end{equation}
converges in $\Rn$ to an $x^*$ fixed point of $F$.
\end{corollaire}

\section{Applications.}\label{app}
\subsection{Solutions of maximal monotone operators.}\label{somm}
In this section, we apply the parallel Jacobi algorithm to the
proximal mapping $F=(I+T)^{-1}$ associated with the maximal
monotone operator $T$. We give first a general result concerning
the maximal monotone operators. Such operators have been studied
extensively because of their role in convex analysis (minimization
of functionals, min-max problems, convex
programming, ...) and certain partial differential equations (Rockafellar[\ref{fellar76a}]).\\
Let $T$ be a multivalued maximal monotone operator defined from
$\Rn$ to $\Rn$. A fundamental problem is to determine an $x^*$ in
$\Rn$ satisfying $0\in Tx^*$ which will be called a solution of
the operator $T$.

\begin{theoreme}\label{jacobioperators}
Let $T$ be a multivalued maximal monotone operator such that
$T^{-1}0\neq\Empty$. Then every parallel Jacobi algorithm
associated with the single-valued mapping $F=(I+T)^{-1}$ converges
in $\Rn$ to an $x^*$ solution of the  problem $0\in Tx$.
\end{theoreme}
\proof
\begin{equation}\label{zerofixe}
\begin{array}{lll}
0\in Tx & \iff & x \in (I+T)x \\
        & \iff & x = (I+T)^{-1}x \\
        & \iff & x =Fx \\
\end{array}
\end{equation}
Thus, the solutions of $T$ are the fixed points of $F$, so the
condition $T^{-1}0\neq\Empty$ implies the existence of a fixed
point $u$ of $\Rn$. It remains to show that $F$ verifies the
condition $(h_3)$ and apply Corollary \ref{jacobi}. Consider $x^i
\in\Rn\ (i=1,2)$ and put $y^i=Fx^i $ then $x^i \in y^i+Ty^i$ or
$x^i-y^i \in Ty^i$. As $T$ is monotone we have
$\ip{(x^1-y^1)-(x^2-y^2),y^1-y^2} \geq 0$ and therefore
$\ip{x^1-x^2,y^1-y^2 } -\norm{y^1-y^2}^2\geq 0$ which implies
$\norm{Fx^1-Fx^2}^2 \le \ip{Fx^1-Fx^2,x^1-x^2}$
\endproof

\subsection{Minimization of functional.}\label{mf}
\begin{corollaire}\label{jacobiminimisation}
Let $f:\Rn \to \R\cup\set{\infty}$ be a lower semicontinuous
convex function which is proper (i.e not identically $+\infty$).
Suppose that the minimization problem $\Min_{\Rn}f(x)$ has a
solution. Then any parallel Jacobi algorithm associated with the
single-valued mapping $F=(I+\partial f)^{-1}$ converges to a
minimizer of $f$ in $\Rn$.
\end{corollaire}
\proof Since in this case the subdifferential $\partial f$ is
maximal monotone. Moreover the minimizers of $f$ are the solutions
of $\partial f$. We then apply Theorem \ref{jacobioperators} to
$\partial f$.
\endproof

\subsection{Saddle point.}\label{sp}
In this paragraph, we apply Theorem \ref{jacobioperators} to
calculate a saddle point of functional $L:\H \to
[-\infty,+\infty]$. Recall that a saddle point of $L$ is an
element $(x^*,y^*)$ of $\H$ satisfying
$$
L(x^*,y)\le L(x^*,y^*)\le L(x,y^*),\ \forall (x,y) \in \H
$$
which is equivalent to
$$
L(x^*,y^*)= \Inf_{x\in\Rn}L(x,y^*)=\Sup_{y\in\Rp}L(x^*,y)
$$
Suppose that $L(x,y)$ is convex lower semicontinuous in $x \in \Rn$ and concave
upper semicontinuous in $y \in \Rp$. Such functionals are called saddle
functions in the terminology of Rockafellar[\ref{fellar70a}]. Let $T_L$ be a
multifunction defined in $\H$ by
$$
(u,v)\in T_L(x,y)\iff
\left\{
\begin{array}{l}
L(x,y')+\ip{y'-y,v} \le L(x,y) \le L(x',y)-\ip{x'-x,u}\\
      \forall (x',y')\in \H
\end{array}
\right.\\
$$

If $L$ is proper and closed in a certain general sense, then $T_L$
is maximal monotone; see
Rockafellar[\ref{fellar70a},\ref{fellar70b}]. In this case the
global saddle points of $L$ (with respect to minimizing in $x$ and
maximizing in $y$) are the elements $(x,y)$ solutions of the
problem $(0,0)\in T_L(x,y)$. That is
\[
(0,0)\in T_L(x^*,y^*) \iff (x^*,y^*)=arg\
\Min_{x\in\Rn}\Max_{y\in\Rp} L(x,y)
\]
We can then apply Theorem \ref{jacobioperators} to the operator $T_L$ so,
\begin{corollaire}
Let $L$ be a proper saddle function from $\H$ into
$[-\infty,+\infty]$ having a saddle point. Then any parallel
Jacobi algorithm associated with the single-valued mapping
$F=(I+T_L)^{-1}$ from $\H$ into $\H$ converges to a saddle point
of $L$.
\end{corollaire}

\subsection{Convex programming.}\label{pc}
We consider now the convex programming problem,
\begin{equation}\label{programmeconvexeprimal}
(P)\left\{
\begin{array}{l}
Min\ f_0(x), \  x\in \Rn \\
f_i(x) \le 0,\ (1\le i\le m)
\end{array}
\right.\\
\end{equation}
where $f_i:\Rn \to\R\ (0\le i\le m)$ is lower semicontinuous
convex functions. This problem can be reduced to an unconstrained
one by mens of the Lagrangian,
$$L(x,y)=f_0(x)+\Sum_{i=1}^my_if_i(x)$$
where $x\in\Rn$ and $y\in(\R_+)^m$. We observe that $L$ is a
saddle function in the sense of [\ref{fellar70a},p. 363], due to
the assumptions of convexity and continuity. The dual problem
associated with $(P)$ is,
\begin{equation}\label{programmeconvexedual}
(D)\left\{
\begin{array}{l}
Max\ \{g_0(y) = \Inf_{x\in\Rn}L(x,y)\} \\
y\in (\R_+)^m
\end{array}
\right.\\
\end{equation}
If $(x^*,y^*)$ is a saddle point of the Lagrangian $L$ then $x^*$ is an
optimal solution of the primal problem $(P)$ and $y^*$ is an optimal
solution of the dual problem $(D)$.\\
Let $\partial L(x,y)$ the subdifferential of $L$ at $(x,y)\in\H$,
be defined as the set of vectors $(u,v)\in\H$ satisfying
$$
\forall (x',y')\in \H,\ L(x,y')-\ip{y'-y,v}\le L(x,y)\le
L(x',y)-\ip{x'-x,u}
$$
(see Luque[\ref{luque84}] and Rockafellar[\ref{fellar70a}]).\\
Then the operator $T_L:(x,y)\to\set{(u,v):(u,-v)\in\partial
L(x,y)}$ is maximal monotone (Rockafellar[\ref{fellar70a}, Cor.
37.5.2]), so we apply Theorem \ref{jacobioperators} to $T_L$.

\begin{corollaire}
Suppose that the convex programming $(P)$ defined by
(\ref{programmeconvexeprimal}) has a solution. Then any parallel
Jacobi algorithm associated with the single-valued mapping
$F=(I+T_L)^{-1}$ from $\H$ to $\H$ converges to a saddle point
$(x^*,y^*)$  of $L$, and so $x^*$ is a solution of the primal
$(P)$ and $y^*$ a solution of the dual $(D)$.
\end{corollaire}

\subsection{Variational inequality.}\label{iv}
A simple formulation of the variational inequality problem is to
find an $x^*\in\Rn$ satisfying
\begin{equation}\label{inegalitevariationnelle}
\ip{Ax^*,x-x^*}\geq 0\ \forall x\in\Rn
\end{equation}
where $A:\Rn \to \Rn$ is a single-valued monotone and maximal
operator\footnote{\ In fact, it's sufficient that $A$ is monotone
and hemicontinuous, i.e verifying $\Lim_{t\to 0^+}\ip{A(x+ty),
h}=\ip{Ax,h}\ \forall x,y,h\in\Rn.$}. Which is equivalent to find
an $x^*\in\Rn$ such that
$$
0 \in Ax^*+N(x^*)
$$
where $N(x)$ is the normal cone to $\Rn$ at $x$ defined by
(see Rockafellar[\ref{fellar70a},\ref{fellar76a}]),
$$
N(x)=\set{y\in \Rn :\ip{y,x-z}\geq 0\ \forall z\in \Rn}
$$
Rockafellar[\ref{fellar76a}] has considered the multifunction
$T$ defined in $\Rn$ by
\begin{equation}\label{operatort}
Tx=Ax+N(x)
\end{equation}
and shown in [\ref{fellar70c}] that $T$ is maximal monotone. The
relation $0\in Tx^*$ is so that reduced to $-Ax^* \in N(x^*)$ or
$\ip{-Ax^*,x^*-z}\geq 0 \ \forall z \in \Rn$ which is the
variational inequality (\ref{inegalitevariationnelle}). Therefore
the solutions of the operator $T$ (defined by (\ref{operatort}))
are exactly the solutions of the variational inequality
(\ref{inegalitevariationnelle}). By using Theorem
\ref{jacobioperators} we can write
\begin{corollaire}\label{jacobiinegalitevraitionnelle}
Let $A:\Rn \to \Rn$ be a single-valued monotone and hemicontinuous
operator such that the problem (\ref{inegalitevariationnelle}) has
a solution, then any parallel Jacobi algorithm associated with the
single-valued mapping $F=(I+T)^{-1}$ where $T$ is defined by
(\ref{operatort}) converges to $x^*$ solution of the problem
(\ref{inegalitevariationnelle}).
\end{corollaire}


\newpage
\noindent Ahmed ADDOU \\
D\'epartement de math\'ematiques et d'informatique\\
Facult\'e des sciences\\
Universit\'e Mohamed Premier\\
60.000 OUJDA\\MAROC\\
e-mail : addou@sciences.univ-oujda.ac.ma\\\\
Abdenasser BENAHMED \\
Lyc\'ee Oued Eddahab Lazaret\\ 60000 Oujda\\ Maroc\\
e-mail : benahmed.univ.oujda@menara.ma

\end{document}